\documentclass[11pt]{article}
\usepackage{amsmath,amssymb,amsthm,latexsym}
\usepackage{indentfirst}
\def\H{\mathbb{H}}
\def\0{{\bf 0}}
\def\a{\alpha}
\def\b{\beta}
\def\B{\mathbb B}
\def\c{{\bf c}}
\def\d{\delta}

\def\g{\gamma}
\def\G{\mathbb G}

\def\l{\lambda}

\def\V{\mathbb V}
\def\P{\mathbb P }
\def\Q{\mathbb Q }
\def\r{\rho}
\def\R{\mathbb R}
\def\s{\sigma}

\def\m{\mu}

\def\o{\omega}
\def\O{\Omega}
\def\e{\varepsilon}

\theoremstyle{definition}

\title{\textbf{Laguerre Geometry of Hypersurfaces in $\R^n$}}
\author {Tongzhu Li
 \and Changping Wang
\thanks {Partially supported by RFDP and No.10125105 of NSFC.}}
\date{}
\begin{document}
\maketitle

\begin{abstract}
Laguerre geometry of surfaces in $\R^3$ is given in the book of
Blaschke [1], and have been studied by E.Musso and L.Nicolodi [5],
[6], [7], B. Palmer [8] and other authors. In this paper we study
Laguerre differential geometry of hypersurfaces in $\R^n$. For any
umbilical free hypersurface $x: M\to\R^n$ with non-zero principal
curvatures we define a Laguerre invariant metric $g$ on $M$ and a
Laguerre invariant self-adjoint operator ${\mathbb S}: TM\to TM$,
and show that $\{g,{\mathbb S}\}$ is a complete Laguerre invariant
system for hypersurfaces in $\R^n$ with $n\ge 4$. We calculate the
Euler-Lagrange equation for the Laguerre volume functional of
Laguerre metric by using Laguerre invariants. Using the Euclidean
space $\R^n$, the Lorentzian space $\R^n_1$ and the degenerate
space $\R^n_0$ we define three Laguerre space forms $U\R^n$,
$U\R^n_1$ and $U\R^n_0$ and define the Laguerre embedding $
U\R^n_1\to U\R^n$ and $U\R^n_0\to U\R^n$, analogue to the Moebius
geometry where we have Moebius space forms $S^n$, $\H^n$ and
$\R^n$ (spaces of constant curvature) and conformal embedding
$\H^n\to S^n$ and $\R^n\to S^n$ (cf. [4], [10]). Using these
Laguerre embedding we can unify the Laguerre geometry of
hypersurfaces in $\R^n$, $\R^n_1$ and $\R^n_0$. As an example we
show that minimal surfaces in $\R^3_1$ or $\R_0^3$ are Laguerre
minimal in $\R^3$.
\end{abstract}

\medskip\noindent
{\bf 2000 Mathematics Subject Classification:} Primary 53A40;
Secondary 53B25.

\medskip\noindent
{\bf Key words and phrases:} Laguerre transformation group,
Laguerre invariants, Laguerre space forms, Laguerre minimal
hypersurfaces.
\par\bigskip\noindent
{\bf {\S} 1. Introduction}
\par\medskip
In study of contact structure in the unit tangent bundle $US^n$
over unit sphere $S^n$ Sophus Lie discovered a interesting finite
dimensional transformation group $LT\G$, which preserves oriented
(n-1)-spheres in $US^n$. This group  $LT\G$ is called Lie sphere
transformation group, which is isomorphic to the group
$O(n+1,2)/\{\pm 1\}$, where $O(n+1,2)$ is the Lorentzian group in
the Lorentzian space $\R^{n+3}_2$. There are two interesting types
of subgroups of  $LT\G$, one is called Moebius group $M\G$,
consisting of all elements of $O(n+1,2)$ which fix a time-like
vector in $\R^{n+3}_2$; another is called Laguerre group $L\G$,
consisting of all elements of $O(n+1,2)$ which fix a light-like
vector in $\R^{n+3}_2$.
\par\medskip
In Laguerre differential geometry we study invariants of
hypersurfaces in Euclidean space $\R^n$ under the Laguerre
transformation group.
\par\medskip
Laguerre geometry of surfaces in $\R^3$ is given in the book of
W.Blaschke [1], and have been studied by E.Musso and L.Nicolodi
[5], [6] and other authors.
\par\medskip
In this paper we study Laguerre differential geometry of
hypersurfaces in $\R^n$. For any umbilical-free hypersurface $x:
M\to\R^n$ with non-zero principal curvatures we define a Laguerre
invariant metric $g$ on $M$ and a Laguerre invariant self-adjoint
operator ${\mathbb S}: TM\to TM$, and show that $\{g,{\mathbb
S}\}$ is a complete Laguerre invariant system for hypersurfaces in
$\R^n$ with $n\ge 4$. We calculate the Euler-Lagrange equation for
the Laguerre volume functional by using Laguerre invariants.
\par\medskip
Using $\R^n$, Lorentzian space $\R^n_1$ and degenerate space
$\R^n_0$ corresponding to the space-like hyperplane, Lorentzian
hyperplane and degenerate hyperplane in $\R^{n+1}_1$ we define
three Laguerre space forms $U\R^n$, $U\R^n_1$ and $U\R^n_0$ as
suitable bundle over $\R^n$, $\R^n_1$, $\R^n_0$ and define the
Laguerre embedding $ U\R^n_1\to U\R^n$ and $U\R^n_0\to U\R^n$,
analogue to the Moebius geometry where we have Moebius space forms
$S^n$, $\H^n$ and $\R^n$ (spaces of constant curvature) and
conformal embedding $\H^n\to S^n$ and $\R^n\to S^n$ (cf. [4],
[10]). Using these Laguerre embedding we can unify the Laguerre
geometries of hypersurfaces in $\R^n$, $\R^n_1$ and $\R^n_0$. As
an example we show that minimal surfaces in $\R^3_1$ or $\R_0^3$
are Laguerre minimal in $\R^3$.
\par\medskip
We organize the paper as follows. In \S 2 we study the geometry of
oriented spheres in $\R^n$. In \S 3 we study Laguerre
transformation group on $U\R^n$. In \S 4 we define Laguerre space
forms and Laguerre embedding. In \S 5 and \S 6 we study Laguerre
invariants for hypersurfaces in $\R^n$ and prove the fundamental
theorem. In \S 7 we calculate Euler-Lagrange equation for volume
function of Laguerre metric. In \S 8 we unify the geometry of
Laguerre hypersurfaces in $\R^n_1$, $\R^n_0$ and $\R^n$.
\par\bigskip\noindent
{\bf {\S} 2. Geometry of oriented spheres in $\R^n$}
\par\medskip
Let $U\R^n$ be the unit tangent bundle over $\R^n$, which is the
hypersurface in $\R^{2n}$:
$$ U\R^n=\{(x,\xi)\mid x\in \R^n, \xi\in S^{n-1}\}=\R^n\times
S^{n-1}\subset \R^{2n}.\leqno (2.1)$$ An oriented sphere in
$U\R^n$ centered at $p$ with radius $r$ is the $(n-1)-$
dimensional submanifold in $U\R^n$ given by
$$ S(p,r)=\{(x,\xi)\in U\R^n\mid x-p=r\xi\}.\leqno (2.2)$$
Geometrically, $S(p,r)$ with $r\not=0$ corresponds to the oriented
sphere in $\R^n$ centered at $p\in \R^n$ with radius $|r|$. If
$r>0$, the unit normal $\xi$ of $S(p,r)$ is outward; if $r<0$, the
unit normal $\xi$ of $S(p,r)$ is inward. If $r=0$, then
$S(p,r)\subset U\R^n$ consists of all unit tangent vector at $p$.
We call $S(p,0)$ the point sphere at $p\in \R^n$. An oriented
hyperplane in $U\R^n$ with constant unit normal $\xi\in S^{n-1}$
and constant $\l\in \R$ is the $(n-1)-$ dimensional submanifold in
$U\R^n$ given by
$$ P(\xi, \l)=\{(x,\xi)\in U\R^n\mid x\cdot\xi=\l\}.\leqno (2.3)$$
Geometrically, it is the hyperplane $\{x\in \R^n\mid
x\cdot\xi=\l\}$ in $\R^n$ with the unit normal $\xi$.
\par\medskip
We denote by $\Sigma$ the set of all oriented spheres and oriented
hyperplanes in $U\R^n$. If $\g_1,\g_2\in\Sigma$ satisfies
$\g_1=\g_2$, or they intersect in a single point $(x,\xi)\in
U\R^n$, we say that $\g_1$ and $\g_2$ are oriented contact.
Geometrically, $\g_1, \g_2$ are oriented contact at $(x,\xi)$ if
and only if they are spheres in $\R^n$ which touch in $x$ with the
same unit normal $\xi$. We note that any point $(x,\xi)\in U\R^n$
determines uniquely a pencil of oriented spheres contact at $x\in
\R^n$ with the common unit normal $\xi$. We note also that there
is a unique point sphere $S(x,0)$ and a unique hyperplane $P(\xi,
x\cdot\xi)$ in this pencil.
\par\medskip
Let $\R^{n+3}_2$ be the space $\R^{n+3}$, equipped with the inner
product
$$ <X,Y>=-X_1Y_1+X_2Y_2+\cdots
+X_{n+2}Y_{n+2}-X_{n+3}Y_{n+3}.\leqno (2.4)$$ Let $C^{n+2}$ be the
light-cone in $\R^{n+3}$ given by
$$ C^{n+2}=\{X\in \R^{n+3}_2\mid <X,X>=0\}.\leqno (2.5)$$
We denote by $\Q^{n+1}$ the quadric in the real projective space
$R\P^{n+2}$, defined by
$$ \Q^{n+1}=\{[X]\in R\P^{n+1}\mid <X,X>=0\}.\leqno (2.6)$$
Then we can assign an oriented sphere $S(p,r)\in \Sigma$ to a
point $[\g]\in\Q^{n+1}$ by
$$ S(p,r)\leftrightarrow [\g],\hskip 5pt \g=(\frac{1}{2}(1+|p|^2-r^2),\frac{1}{2}(1-|p|^2+r^2),p,
-r)\leqno (2.7)$$ and assign an oriented hyperplane
$P(\xi,\l)\in\Sigma$ to a point in $[\g]\in\Q^{n+1}$ by
$$P(\xi,\l)\leftrightarrow [\g],\hskip 5pt \g=(\l,-\l,\xi,1).\leqno (2.8)$$
We call $[\g]\in \Q^{n+1}$ the coordinate of the oriented sphere
$S(p,r)$ or $P(\xi,\l)$.
\par\medskip
For any $\g\in \Sigma$ we will denote by $[\g]\in\Q^{n+1}$ its
coordinate given in (2.7) and (2.8). It is easy to verify that the
corresponding $ \g\in\Sigma\to [\g]\in\Q^{n+1}$ defines a
bijection from $\Sigma$ to $\Q^{n+1}\backslash\{[\wp]\}$, where
$$\wp=(1,-1,\0,0),\hskip 5pt \0\in \R^n.\leqno (2.9)$$
Geometrically, the point
$$[\wp]=\lim_{|p|\to\infty}\frac{1}{|p|^2}[(\frac{1}{2}(1+|p|^2),\frac{1}{2}(1-|p|^2),p,
0)]$$ in $\Q^{n+1}$ is the coordinate of the point sphere at
$\infty$ of $\R^n$. Using (2.7) and (2.8) we can easily verify
that $\g_1,\g_2\in\Sigma$ are oriented contact if and only if
their sphere coordinates $[\g_1]$ and $[\g_2]$satisfy
$$ <\g_1,\g_2>=0.\leqno (2.10)$$
From (2.7), (2.8) and (2.9) we know that $[\g]\in
\Q^{n+1}\backslash\{[\wp]\}$ is an oriented sphere $S(p,r)$ in
$U\R^n$ if and only if $<\g,\wp>\not=0$ and that $[\g]\in
\Q^{n+1}\backslash\{[\wp]\}$ is a hyperplane $P(\xi,\l)$ in
$U\R^n$ if and only if $<\g,\wp>=0$.
\par\medskip
A point $(x,\xi)\in U\R^n$ determines a unique pencil of oriented
spheres contact at $x\in \R^n$ with the common unit normal $\xi$,
and the point sphere $\g_1=S(x,0)$ and the oriented hyperplane
$\g_2=P(\xi,x\cdot\xi)$ in the pencil have coordinate $[\g_1]$ and
$[\g_2]$, where
$$ \g_1=(\frac{1}{2}(1+|x|^2),\frac{1}{2}(1-|x|^2),x,
0),\hskip 5pt \g_2=(x\cdot\xi,-x\cdot\xi,\xi,1).\leqno (2.11)$$
Then any oriented sphere $[\g]$ in the pencil can be written as
$$ [\g]=[\l \g_1+\m \g_2]\in  \Q^{n+1}\backslash\{[\wp]\}\leqno
(2.12)$$ for some $(\l,\m)\in \R^2\backslash\{0\}$. Thus a point
$(x,\xi)\in U\R^n$ determines a unique projective line
$$ \{[\l \g_1+\m \g_2]\mid (\l,\m)\in \R^2\backslash\{0\}\}$$
lying in $\Q^{n+1}\backslash\{[\wp]\}$.
\par\medskip
Let $\Lambda^{2n-1}$ be the set consisting of all projective lines
lying in $\Q^{n+1}\backslash\{[\wp]\}$. Then the mapping $L:
U\R^n\to\Lambda^{2n-1}$ defined by
$$ L((x,\xi))=\{[\l{\g_1}+\m {\g_2}]\in \Q^{n+1}\backslash\{[\wp]\}
 \mid (\l,\m)\in \R^2\backslash\{0\}\}\leqno (2.13)$$
is a diffeomorphism, called Lie diffeomorphism.
\par\bigskip\noindent
{\bf {\S} 3. Laguerre transformation group on $U\R^n$}
\par\medskip
Let $\G$ be the subgroup of Lorentzian group $O(n+1,2)$ on
$\R^{n+3}_2$ given by
$$ L\G=\{T\in O(n+1,2)\mid \wp T=\wp\},\leqno (3.1)$$
where $\wp$ is the light-like vector in $\R^{n+3}_2$ defined by
(2.9). Then any $T\in L\G$ induces a transformation on $\Q^{n+1}$
defined by
$$ T([X])=[XT],\hskip 5pt X\in \Q^{n+1}.\leqno (3.2)$$
We call both $T\in L\G$ and $T:\Q^{n+1}\to\Q^{n+1}$ Laguerre
transformation.
\par\medskip
Let $\g_1,\g_2\in\Sigma$ be two different oriented contact spheres
or hyperplanes. Then $[\g_1]$ and $[\g_2]$ define a projective
line lying in $\Q^{n+1}\backslash\{[\wp]\}$ by
$$ span\{[\g_1], [\g_2]\}=\{[\l{\g_1}+\m {\g_2}]
 \mid (\l,\m)\in \R^2\backslash\{0\}\}\in \Lambda^{2n-1}.$$
Then any $T\in L\G$ defines a transformation $T:\Lambda^{2n-1}\to
\Lambda^{2n-1}$ by
$$ T(span\{[\g_1],[\g_2]\})=span\{[\g_1T],[\g_2T]\}. $$
Let $L: U\R^n\to  \Lambda^{2n-1}$ be the Lie diffeomorphism. Then
any $T\in L\G$ induces a transformation
$$ \s=L^{-1}\circ T\circ L: U\R^n\to U\R^n, $$
called a Laguerre transformation on $U\R^n$. Thus the Laguerre
transformation group on $U\R^n$ is given by
$$ L\G=\{\s: U\R^n\to U\R^n\mid \s=L^{-1}\circ T\circ L,T\in O(n+1,2),
\wp T=\wp\}.$$ The dimension of $L\G$ is $(n+2)(n+1)/2$.
\par\medskip
Let $T\in L\G$ be a Laguerre transformation. Then we have $\wp
T=\wp$ and $T: \Q^{n+1}\backslash\{[\wp]\}\to
\Q^{n+1}\backslash\{[\wp]\}$. Since any oriented sphere $\g\in
U\R^n$ determines a point $[\g]\in \Q^{n+1}\backslash\{[\wp]$ such
that $<\g,\wp>\not=0$, then $<\g T,\wp>=<\g T,\wp
T>=<\g,\wp>\not=0$, we know that $T([\g])$ is also an oriented
sphere. Since any oriented hyperplane $\g\in U\R^n$ determines a
point $[\g]\in \Q^{n+1}\backslash\{[\wp]$ such that $<\g,\wp>=0$,
then $<\g T,\wp>=<\g T,\wp T>=<\g,\wp>=0$, we know that $T([\g])$
is also an oriented hyperplane. Thus any Laguerre transformation
$\s: U\R^n\to U\R^n$ takes oriented spheres to oriented spheres,
takes oriented hyperplanes to oriented hyperplanes.
\par\medskip\noindent
{\bf Example 3.1.} Any isometry in $\R^n$ given by
$$ \s (x)= xA+a,\hskip 5pt A\in O(n), a\in \R^n$$
induces an isometry transformation $\s : U\R^n\to U\R^n$ defined
by
$$ \s ((x,\xi))=(xA+a, \xi A).\leqno (3.3)$$
It is easy to check that $\s$ is a Laguerre transformation on
$U\R^n$ and that
$$ T(\s)=L\circ\s\circ
L^{-1}=\left(\begin{matrix}1+|a|^2/2&-|a|^2/2&a&0\\|a|^2/2&1-|a|^2/2&a&0\\Aa'&
-Aa'&A&0\\0&0&0&1\end{matrix}\right)\in L\G,\leqno (3.4)$$ where
$a'$ is the transport of the vector $a\in\R^n$.
\par\medskip\noindent
{\bf Example 3.2.} The 1-parametric parabolic transformations (or
parallel transformations) defined by
$$ \phi_t(x,\xi)=(x+t\xi,\xi),\hskip 5pt t\in \R.\leqno (3.5)$$
are Laguerre transformations in $U\R^n$. It is easy to check that
$$ T(\phi_t)=L\circ\phi_t\circ L^{-1}=\left(\begin{matrix}1-t^2/2&t^2/2&0&-t\\-t^2/2&1+t^2/2&0&-t\\0&0&I_n&0\\
t&-t&0&1\end{matrix}\right)\in L\G.\leqno (3.6)$$ Since
$\phi_s\circ\phi_t=\phi_{s+t}$, we call $\phi_t$ a parabolic flow
in $U\R^n$.
\par\medskip\noindent
{\bf Example 3.3.} The third example of Laguerre transformations
in $U\R^n$ is the following 1-parametric hyperbolic
transformations. For any $(x,\xi)\in U\R^n$ we write
$$ x=(x_0,x_1)\in \R^{n-1}\times\R,\hskip 5pt
\xi=(\xi_0,\xi_1)\in \R^{n-1}\times\R,$$ then a hyperbolic
transformation $$\psi_t (x,\xi)=(x(t),\xi (t))\in U\R^n,\hskip 5pt
t\in\R$$ is defined by
$$ x(t)=\left( x_0-\frac{\sinh t
x_1}{\sinh t\xi_1+\cosh t}\xi_0,\frac{x_1}{\sinh t\xi_1+\cosh
t}\right);\leqno (3.7)$$
$$ \xi (t)=\left(\frac{1}{\sinh t\xi_1+\cosh t}\xi_0 ,\frac{\cosh t\xi_1+\sinh t}{\sinh t\xi_1+\cosh t}
\right).\leqno (3.8)$$ It is easy to check that
$$ T(\psi_t)=L\circ\psi_t\circ
L^{-1}=\left(\begin{matrix}I_{n+1}&0&0\\0&\cosh t&\sinh t
\\0&\sinh t&\cosh t
\end{matrix}\right)\in L\G.\leqno (3.9)$$ Since
$\psi_s\circ\psi_t=\psi_{s+t}$, we call $\psi_t$ a hyperbolic flow
in $U\R^n$.
\par\medskip
Let $\{e_1,e_2,\cdots, e_{n+3}\}$ be the standard basis for
$\R^{n+3}_2$, $e_i=(0,\cdots,0,1,0,\cdots,0)$. For any $T\in L\G$
we have
$$ \wp T=\wp, \hskip 5pt <e_iT,\wp>=<e_iT,\wp T>=<e_i,\wp>, 1\le
i\le n+3.$$ Using these information and the fact that $T\in
O(n+1,2)$ we can write
$$
T=\left(\begin{matrix}1+|a|^2/2-\r^2/2&-|a|^2/2+\r^2/2&a&\r\\|a|^2/2-\r^2/2&1-|a|^2/2+\r^2/2&a&\r
\\Aa'-\r u&-Aa'+\r u&A&u\\va'-\r w&-va'+\r
w&v&w\end{matrix}\right)\leqno (3.10)$$ for some $$
\left(\begin{matrix}A&u\\v&w\end{matrix}\right)\in O(n,1),\hskip
5pt (a,\r)\in \R^{n+1}, w\in\R.\leqno (3.11)$$ It is easy to check
that
$$ T\to
\left(\begin{matrix}A&u&0\\v&w&0\\a&\r&1\end{matrix}\right)\leqno
(3.12)$$ is a isomorphism from $L\G$ to the Lorentzian
transformation group in $\R^{n+1}_1$.
\par\medskip
Now let $\g_1=S(p,r)$, $\g_2=S(p^*,r^*)$ be oriented spheres in
$\R^n$. Let $T$ be a Laguerre transformation given by (3.10).
Since
$$ \g_1=(\frac{1}{2}(1+|p|^2-r^2),\frac{1}{2}(1-|p|^2+r^2),p,
-r),$$$$
\g_2=(\frac{1}{2}(1+|p^*|^2-r^{*2}),\frac{1}{2}(1-|p^*|^2+r^{*2}),p^*,
-r^*),$$ then the oriented spheres $\g_1T=S({\tilde p},{\tilde
r})$ and $\g_2T=S({\tilde p^*},{\tilde r^*})$ are given by
$$ ({\tilde p},-{\tilde
r})=(pA-rv+a,pu+rw+\r), \hskip 5pt ({\tilde p^*},-{\tilde
r^*})=(p^*A-r^*v+a,p^*u+r^*w+\r).$$ Thus we have
$$ ({\tilde p^*}-{\tilde p}, -{\tilde r^*}+{\tilde
r})=(p^*-p,-r^*+r)\left(\begin{matrix}A&u\\v&w\end{matrix}\right).\leqno
(3.13)$$ It follows that
$$ F=|p^*-p|^2-(r^*-r)^2\leqno (3.14)$$
is a Laguerre invariant. Geometrically,  if one sphere is not
contained in another, then $F$ is exactly the square length of the
common tangent segment of the two spheres $S(p,r)$ and
$S(p^*,r^*)$.
\par\medskip
{\bf Theorem 3.1}{\hskip 3pt}{\it For any $T\in O(n+1,2)$ with
$\wp T=T$ there exist two isometries $\s_1$, $\s_2$ on $U\R^n$ and
constants $s,t\in\R, \e=\pm 1$ such that
$$ T=\e T(\s_2)T(\psi_t)T(\phi_s)T(\s_1).\leqno (3.15)$$ }
\par\smallskip\noindent
{\bf Proof.}{\hskip 3pt} For any $T\in O(n+1,2)$ with $\wp T=T$ we
can write $T$ as in (3.10). From (3.11) we get $w^2=1+|v|^2$. By
changing $T$ to $-T$, if necessary, we may assume that $w>0$. Then
we have
$$ e_{n+3}T=((va'-\r w),-(va'-\r w),v,w),\hskip 5pt w=\sqrt{|v|^2+1}.$$
We can find $s,t\in\R$ and $A_1\in O(n)$
$$ s=w^{-1}(va'-\r w),\hskip 5pt w=\cosh t,\hskip 5pt vA_1=(0,\cdots,0,\sinh t).$$
We denote by $\s_1^{-1}$ the isometry
$\s_1^{-1}((x,\xi))=(xA_1,\xi A_1)$ on $U\R^n$, then by (3.4),
(3.6) and (3.9) we have
$$ e_{n+3}TT(\s_1^{-1})T(\phi_{-s})T(\psi_{-t})=e_{n+3}.$$
We define $$T^*= TT(\s_1^{-1})T(\phi_{-s})T(\psi_{-t}).$$ Since
$T^*$ satisfies
$$T^*\in O(n+1,2),\hskip 5pt \wp
T^*=\wp,\hskip 5pt e_{n+3}T^*=e_{n+3},$$ $T^*$ takes the form
(3.10) with $e_{n+3}=(0,\cdots,0,1)$ as its last line. Thus $T^*$
takes the form (3.4) for some isometry $\s_2$. Thus we get (3.15)
and complete the proof of Theorem 2.1.
\par\medskip
{\bf Corollary}{\hskip 3pt}{\it Any Laguerre transformation in
$U\R^n$ is generated by the isometries, the parallel
transformations and the hyperbolic transformations.}
\par\bigskip\noindent
{\bf {\S} 4. Laguerre space forms and Laguerre embeddings}
\par\medskip
In Moebius geometry we have three standard spaces $S^n$, $\R^n$
and $\H^n$ and conformal embeddings $\R^n\to S^n$ and $\H^n\to
S^n$. Similarly, we introduce in this section three Laguerre space
forms $U\R^n$, $U\R^n_1$ and $U\R^n_0$ and the Laguerre embeddings
$\s: U\R^n_1\to U\R^n$ and $\tau: U\R^n_0\to U\R^n$.
\par\medskip
Let $\R^n_1$ be the Lorentzian space with inner product
$$ <X,Y>=X_1Y_1+\cdots+X_{n-1}Y_{n-1}-X_nY_n.\leqno (4.1)$$
Let $U\R^n_1$ be the unit bundle of $\R^n_1$ defined by
$$ U\R_1^n=\{(x,\xi)\mid x\in \R_1^n, <\xi,\xi>=-1\}.\leqno (4.2)$$
\par\medskip
An oriented sphere (hyperboloid) $H(p,r)$ centered at $p$ in
$\R_1^n$ with radius $r$ can be embedded in $U\R_1^n$ as the
(n-1)-dimensional submanifold
$$ H(p,r)=\{(x,\xi)\in U\R_1^n\mid x-p=r\xi\}.\leqno (4.3)$$
Here $r$ is a real number. If $r=0$, then $H(p,r)$ consists all
unit time-like vectors at $p$, called "point sphere" at $p$. We
assign $\g=H(p,r)$ a vector $[\g]\in \Q^{n+1}$ by
$$
\g=(\frac{1}{2}(1+<p,p>+r^2),\frac{1}{2}(1-<p,p>-r^2),-r,p).\leqno
(4.4)$$
\par\medskip
An oriented space-like hyperplane $P(\xi,\l)$ in $\R^n_1$ with
unit normal $\xi$ can be embedded in $U\R_1^n$ as the
(n-1)-dimensional submanifold
$$ P(\xi,\l)=\{(x,\xi)\in U\R_1^n\mid <x,\xi>=\l\}.\leqno (4.5)$$
We assign $\g=P(\xi,\l)$ the vector $[\g]\in\Q^{n+1}$ by
$$\g=(\l,-\l,1,\xi).\leqno (4.6)$$
\par\medskip
Two oriented spheres (or hyperplanes) $\g_1$ and $\g_2$ are
oriented contact in $\R^n_1$ if and only if their corresponding
vectors $\g_1,\g_2\in\Q^{n+1}$ satisfy $<\g_1,\g_2>=0$. Any point
$(x,\xi)\in U\R^n_1$ determines a pencil of spheres (hyperplanes)
in $\R^n_1$ oriented contacted at $x$ with common normal vector
$\xi$, which corresponds to a projective line by the Lie
diffeomorphism $L_1: U\R^n_1\to \Lambda^{2n-1}$ given by
$$ L_1(x,\xi)=\{[\l\g_1+\m\g_2]\mid (\l,\mu)\in \R^2\backslash\{0\}\},\leqno (4.7)$$
where$$
\g_1=(\frac{1}{2}(1+<x,x>),\frac{1}{2}(1-<x,x>),0,x),\leqno
(4.8)$$$$ \g_2=(<x,\xi>,-<x,\xi>,1,\xi).\leqno (4.9)$$ Let $L:
U\R^n\to \Lambda^{2n-1}$ be the Lie diffeomorphism defined by
(2.13). It is easy to check that $\s=L^{-1}\circ L_1: U\R^n_1\to
U\R^n$ given by
$$ \s (x,\xi)=(x',\xi')\in U\R^n; \leqno (4.10)$$
where $(x,\xi)\in U\R^n_1$ with $x=(x_0,x_1)\in\R^{n-1}\times\R$,
$\xi=(\xi_0,\xi_1)\in\R^{n-1}\times\R$ and
$$x'=(-\frac{x_1}{\xi_1},x_0-\frac{x_1}{\xi_1}\xi_0), \hskip
5pt\xi'=(\frac{1}{\xi_1},\frac{\xi_0}{\xi_1}).\leqno (4.11)$$ It
is straightforward to verify that $\s$ takes the hyperplane
$P(\xi,\l)$ in $U\R^n_1$ to the hyperplane $P(\xi',\l/\xi_1)$ in
$U\R^n$, takes the oriented sphere $H(p,r)$ in $U\R^n_1$ into the
oriented sphere $S(p',r')$ in $U\R^n$, where $p=(p_0,p_1)$,
$p'=(-r,p_0)$ and $r'=-p_1$. Thus $\s:U\R^n_1\to U\R^n$ is a
Laguerre embedding.
\par\medskip
Let $\R^{n+1}_1$ be the Lorentzian space with inner product
$$ <X,Y>=X_1Y_1+\cdots+X_nY_n-X_{n+1}Y_{n+1}.\leqno (4.12)$$
Let $\nu=(1,\0,1)$ be the light-like vector in $\R^{n+1}_1$ with
$\0\in \R^{n-1}$. Let $\R^n_0$ be the degenerate hyperplane in
$\R^{n+1}_1$ defined by
$$ \R^n_0=\{X\in \R^{n+1}_1\mid <X,\nu>=0\}. \leqno (4.13)$$
We define
$$ U\R^n_0=\{(x,\xi)\in \R^{n+1}_1\times\R^{n+1}_1\mid  <x,\nu>=0, <\xi,\xi>=0, <\xi,\nu>=1\}.\leqno (4.14)$$
\par\medskip
An oriented sphere $C(p)$ in $\R_0^n$ with $p\in\R^{n+1}_1$ is the
(n-1)-submanifold in $U\R^n_0$ given by
$$ C(p)=\{(x,\xi)\in U\R^n_0\mid x-p=-<p,\nu>\xi\}.\leqno (4.15)$$
Geometrically, $C(p)$ ($<p,\nu>\not=0$) is a paraboloid in
$\R^n_0$ as the intersection of the light-cone $<X-p,X-p>=0$ in
$\R^{n+1}_1$ with the degenerate hyperplane $<X,\nu>=0$. The
paraboloid $C(p)$ is centered at $p^*=p+(r,\0,0)\in \R^n_0$
($r=-<p,\nu>$) with the symmetric axe $\ell=\{p^*+t\nu\mid
t\in\R\}$. If $r=0$, then $C(p)$ consists of all $(p,\xi)\in
U\R^n_0$ with $\xi$ lying on the paraboloid $\{\xi\in
\R^{n+1}_1\mid <\xi,\xi>=0, <\xi,\nu>=1\}$ in $ \R^{n+1}_1$. We
assign $\g=C(p)$ a vector $[\g]\in \Q^{n+1}$ by
$$
\g=(\frac{1}{2}(1+<p,p>),\frac{1}{2}(1-<p,p>),p).\leqno (4.16)$$
\par\medskip
An oriented space-like hyperplane $P(\xi,\l)$ in $\R^n_0$ with
unit normal $\xi$ can be embedded in $U\R^n_0$ as the
(n-1)-submanifold
$$ P(\xi,\l)=\{(x,\xi)\in U\R^n_0\mid <x,\xi>=\l\}.\leqno (4.17)$$
We assign $\g=P(\xi,\l)$ the vector $[\g]\in\Q^{n+1}$ by
$$\g=(\l,-\l,\xi).\leqno (4.18)$$
\par\medskip
Two oriented spheres (or hyperplanes) $\g_1$ and $\g_2$ are
oriented contact in $\R^n_0$ if and only if their corresponding
vectors $\g_1,\g_2\in\Q^{n+1}$ satisfy $<\g_1,\g_2>=0$. Any point
$(x,\xi)\in U\R^n_0$ determines a pencil of spheres (hyperplanes)
in $\R^n_0$ oriented contacted at $x$ with common normal vector
$\xi$, which corresponds to a projective line by the Lie
diffeomorphism $L_0: U\R^n_0\to \Lambda^{2n-1}$
$$ L_0(x,\xi)=\{[\l\g_1+\m\g_2]\mid (\l,\m)\in\R^2\backslash\{0\}\},\leqno (4.19)$$
where$$ \g_1=(\frac{1}{2}(1+<x,x>),\frac{1}{2}(1-<x,x>),x),\leqno
(4.20)$$$$ \g_2=(<x,\xi>,-<x,\xi>,\xi).\leqno (4.21)$$ Let $L:
U\R^n\to \Lambda^{2n-1}$ be the Lie diffeomorphism defined by
(2.13). It is easy to check that $\tau=L^{-1}\circ L_0: U\R^n_1\to
U\R^n$ given by
$$ \tau (x,\xi)=(x',\xi')\in U\R^n; \leqno (4.22)$$
where $x=(x_1,x_0,x_1)\in\R\times\R^{n-1}\times\R$,
$\xi=(\xi_1+1,\xi_0,\xi_1)\in\R\times\R^{n-1}\times \R$ and
$$x'=(-\frac{x_1}{\xi_1},x_0-\frac{x_1}{\xi_1}\xi_0), \hskip
5pt\xi'=(1+\frac{1}{\xi_1},\frac{\xi_0}{\xi_1}).\leqno (4.23)$$ It
is straightforward to verify that $\tau$ takes the hyperplane
$P(\xi,\l)$ in $U\R^n_0$ to the hyperplane $P(\xi',\l/\xi_1)$ in
$U\R^n$, takes the oriented sphere $C(p)$ in $U\R^n_0$ into the
oriented sphere $S(p',r')$ in $U\R^n$, where $p=(p_1-r,p_0,p_1)$,
$p'=(p_1-r,p_0)$, $r=-<p,\nu>$ and $r'=-p_1$. Thus
$\tau:U\R^n_0\to U\R^n$ is a Laguerre embedding.

\par\bigskip\noindent
{\bf {\S} 5. Laguerre hypersurfaces in $U\R^n$}
\par\medskip
Let $x: U\R^n\to \R^n,\hskip 5pt \xi: U\R^n\to S^{n-1}\subset
\R^n$ be the standard projections $(x,\xi)\to x$ and $(x,\xi)\to
\xi$, respectively. Then there is a standard contact form $\o$ in
$U\R^n$ defined by
$$ \o=dx\cdot\xi.\leqno (5.1)$$
It is easy to verify that $\o\wedge d\o^{n-1}\not=0$, which is (up
to a non-zero constant) the volume form for the embedding of
$U\R^n=\R^n\times S^{n-1}$ in $\R^{2n}$.
\par\medskip
Let $(x,\xi): U\R^n\to \R^n\times S^{n-1}\subset \R^{2n}$ be the
standard embedding. We define $\g_1, \g_2: U\R^n\to \R^{n+3}_2$ by
(2.11). Let $T\in L\G$ be a Laguerre transformation and
$$({\tilde x},{\tilde \xi})=\phi((x,\xi)),\hskip 5pt \phi=L^{-1}\circ T\circ L: U\R^n\to U\R^n. $$
We denote by $a, b$ the last coordinate of $\g_1T$ and $\g_2T$,
respectively. Then by (2.11) and (3.10) we can write
$$ {\tilde \g}_1=(\frac{1}{2}(1+|{\tilde x}|^2),\frac{1}{2}(1-|{\tilde x}|^2),{\tilde x},
0)={\g}_1T-\frac{a}{b}{\g}_2T,\leqno (5.2)$$ $$
{\tilde\g}_2=({\tilde x}\cdot{\tilde\xi},-{\tilde
x}\cdot{\tilde\xi},{\tilde\xi},1)=\frac{1}{b} {\g}_2T.\leqno
(5.3)$$ It follows that
$$ d{\tilde x}\cdot{\tilde\xi}=<d{\tilde\g}_1,{\tilde\g}_2>=<d({\g}_1T-\frac{a}{b}{\g}_2T),\frac{1}{b}
{\g}_2T>=\frac{1}{b}<d\g_1,\g_2>=\frac{1}{b}\,dx\cdot\xi.\leqno
(5.4)$$
$$ d{\tilde \xi}\cdot
d{\tilde\xi}=<d{\tilde\g}_2,d{\tilde\g}_2>=\frac{1}{b^2}<d\g_2,d\g_2>=\frac{1}{b^2}\,d\xi\cdot
d\xi.\leqno (5.5)$$
\par\medskip
We call $f=(x,\xi): M^{n-1}\to U\R^n$ a Laguerre hypersurface, if
$\xi: M^{n-1}\to \R^n$ is a immersion and $ f^*\o=dx\cdot\xi=0$.
It follows from (5.4) and (5.5) that  any Laguerre transformation
takes Laguerre hypersurfaces in $U\R^n$ to Laguerre hypersurfaces
in $U\R^n$. By (2.2) and (2.3) we know that oriented spheres and
hyperplanes are simplest Laguerre hypersurfaces in $U\R^n$.
\par\medskip
Let $x: M^{n-1}\to\R^n$ be an oriented hypersurface in $\R^n$ with
non-zero principal curvatures. Then the unit normal $\xi: M\to
\R^n$ is a immersion. Thus $x$ induces uniquely a Laguerre
hypersurface $f=(x,\xi): M^{n-1}\to U\R^n$. We note that for a
Laguerre hypersurface $f=(x,\xi): M^{n-1}\to U\R^n$, $x:
M^{n-1}\to\R^n$ may not be an immersion. By a theorem of U.
Pinkall [9] we know that the parallel transformation
$f_t=(x+t\xi,\xi)$ of $f$ is an immersion at any given point $p\in
M^{n-1}$ for almost all $t\in\R$. In this sense we may assume that
$x: M^{n-1}\to \R^n$ is an immersion.
\par\medskip
Let $x,{\tilde x}: M^{n-1}\to \R^n$ be two oriented hypersurfaces
with non-zero principal curvatures. We say $x, {\tilde x}$ are
Laguerre equivalent, if the corresponding Laguerre hypersurfaces
$f=(x,\xi), {\tilde f}=({\tilde x},{\tilde \xi}): M^{n-1}\to
U\R^n$ are differ only by a Laguerre transformation $\phi:
U\R^n\to U\R^n$, i.e., ${\tilde f}=\phi\circ f$. In Laguerre
differential geometry we study properties of Laguerre
hypersurfaces in $U\R^n$ which are invariant under the Laguerre
transformation group in $U\R^n$.
\par\medskip
Let $x: M\to \R^n$ be oriented hypersurface with unit normal
$\xi$. We define
$$ [y]: M\to Q^{n+1},\hskip 5pt y=(x\cdot\xi, -x\cdot\xi,\xi,1).\leqno (5.6)$$
\par\medskip
{\bf Theorem 5.1}{\hskip 3pt}{\it Let $x,{ x^*}: M\to \R^n$ be two
oriented hypersurfaces with non-zero principal curvatures. Then
$x$ and ${ x^*}$ are Laguerre equivalent if and only if there
exists $T\in L\G$ such that $[{ y^*}]=[yT]$.}
\par\smallskip\noindent
{\bf Proof.}{\hskip 3pt} Let $\xi$ and $\xi^*$ be the unit normal
of $x$ and $x^*$, respectively. If there is a Laguerre
transformation $\phi=L^{-1}\circ T\circ L\in L\G$ such that
$(x^*,\xi^*)=\phi\circ (x,\xi)$, then by (5.3) we obtain
$[y^*]=[yT]$. Conversely, if $[y^*]=[yT]$ for some $T\in L\G$, we
define $({\tilde x},{\tilde \xi})=\phi\circ (x,\xi)$ with
$\phi=L^{-1}\circ T\circ L$. Then by (5.3) we have $[{\tilde
y}]=[yT]=[y^*]$. It follows that
$$ ({\tilde x}\cdot{\tilde \xi},-{\tilde x}\cdot{\tilde
\xi},{\tilde\xi},1)=(x^*\cdot\xi^*,-x^*\cdot\xi^*,\xi^*,1).\leqno
(5.7)$$ Let $\{e_i\}$ be a local basis for $TM$. Since $\xi^*:
M\to\R^n$ is an immersion, we know that
$\{e_1(\xi^*),\cdots,e_{n-1}(\xi^*),\xi^*\}$ is a basis for
$\R^n$. From (5.7) and the facts that
$$ \xi^*={\tilde\xi},\hskip 5pt (x^*-{\tilde x})\cdot
\xi^*=0,\hskip 5pt (x^*-{\tilde x})\cdot d\xi^*=d((x^*-{\tilde
x})\cdot \xi^*)=0, $$ we get $x^*={\tilde x}$. Thus we have
$(x^*,\xi^*)=\phi\circ (x,\xi)$, which implies that $x$ and ${
x^*}$ are Laguerre equivalent. We complete the proof of Theorem
5.1.
\par\medskip
Since by (5.6) we have $<dy,dy>=d\xi\cdot d\xi$, which is exactly
the third fundamental from of $x$. It follows from Theorem 5.1
that
\par\medskip
{\bf Corollary}{\hskip 3pt}{\it The conformal class of the third
fundamental form of a hypersurface $x: M\to \R^n$ is a Laguerre
invariant.}
\par\medskip
Let $x: M\to \R^n$ be a oriented hypersurface with non-zero
principal curvatures. Let $III=<dy,dy>$ be the third fundamental
form of $x$. For any orthonormal basis $\{E_1,E_2,\cdots,
E_{n-1}\}$ with respect to $III$ we define
$$ \V=span\{y,\Delta y, E_1(y),E_2(y),\cdots, E_{n-1}(y)\},\leqno (5.8)$$
where $\Delta$ is the Laplacian operator with respect to
$III=<dy,dy>$. Then we have
$$ <y, E_i(y)>=<\Delta y, E_i(y)>=0, \hskip 5pt <y,\Delta y>=-(n-1),\hskip 5pt
<E_i(y),E_j(y)>=\d_{ij}.\leqno (5.9)$$ Thus at each point $\V$ is
a $(n+1)-$dimensional non-degenerate subspace in $\R^{n+3}_2$ of
type $(-,+,\cdots,+)$. Let
$$ \R^{n+3}_2=\V\oplus \V^{\perp}=span\{y,\Delta y, E_1(y),E_2(y),\cdots,
E_{n-1}(y)\}\oplus \V^{\perp}\leqno (5.10)$$ be the orthogonal
decomposition of $\R^{n+3}_2$. Then $\V^{\perp}$ is a
2-dimensional non-degenerate subspace of $\R^{n+3}_2$ of type
$(-,+)$.
\par\medskip
Let $\{e_1,e_2,\cdots, e_{n-1}\}$ be the orthonormal basis for
$TM$ with respect to $dx\cdot dx$, consisting of unit principal
vectors. We write the structure equation of $x: M\to\R^n$ by
$$ e_j(e_i(x))=\sum_{k}\Gamma^k_{ij}e_k(x)+k_i\d_{ij}\xi;\hskip 5pt
e_i(\xi)=-k_ie_i(x),\hskip 5pt 1\le i,j,k\le n-1,\leqno (5.11)$$
where $k_i\not=0$ is the principal curvature corresponding to
$e_i$. Let
$$ r_i=\frac{1}{k_i}, \hskip 5pt
r=\frac{r_1+r_2+\cdots+r_{n-1}}{n-1}\leqno (5.12)$$ the curvature
radius and mean curvature radius of $x$. Then the mean curvature
sphere $S(x+r\xi,r)$ of $x$ in $\R^n$ has the sphere coordinate
$[\eta]$, where
$$ \eta =(\frac{1}{2}(1+|x|^2),\frac{1}{2}(1-|x|^2),
x,0)+r\,(x\cdot\xi, -x\cdot\xi,\xi,1).\leqno (5.13)$$ We define
$E_i=r_ie_i, 1\le i\le (n-1)$, then $\{E_1, E_2,\cdots, E_{n-1}\}$
is an orthonormal basis for $III=<dy,dy>$. From (5.11) we get
$$ E_i(y)=-(x\cdot e_i(x),-x\cdot e_i(x),e_i(x), 0).\leqno (5.14)$$
It follows from (5.6), (5.13) and (5.14) that
$$ <y,\eta>=0,\hskip 5pt <E_i(y),\eta>=0,\hskip 5pt
E_i(\eta)=(r-r_i)E_i(y)+E_i(r)y. \leqno (5.15)$$ Moreover, from
(5.15) we get
$$ <\Delta y,
\eta>=\sum_i<E_iE_i(y),\eta>=-\sum_i<E_i(y),E_i(\eta)>=-\sum_i(r-r_i)=0.$$
Thus we know that $\eta\in\V^{\perp}$. Let $\wp=(1,-1,\0,0)\in
\R^{n+3}_2$ be the vector defined by (2.9). Since $<y,\wp>=0$, by
(5.13) we have
$$ \V^{\perp}=span\{\eta,\wp\},\hskip 5pt
<\eta,\eta>=<\wp,\wp>=0,\hskip 5pt <\eta,\wp>=-1.\leqno (5.16)$$
We call $\eta: M\to C^{n+2}\subset \R^{n+3}$ defined by (5.13) the
Laguerre Gauss map of $x$.
\par\medskip
It is clear that $\V$, $\V^{\perp}$ and $\eta$ are Laguerre
invariant: if $x$ is Laguerre equivalent to ${\tilde x}$ by $T\in
L\G$, then we have
$${\tilde \V}=\V T,\hskip 5pt {\tilde
\V}^{\perp}=\V^{\perp} T,\hskip 5pt {\tilde\eta}=\eta T.\leqno
(5.17)$$
\par\medskip
Now let $x,{\tilde x}: M\to \R^n$ are Laguerre equivalent by $T\in
L\G$. Then by (5.3) and (5.17) we have
$$ {\tilde y}=\frac{1}{b}yT,\hskip 5pt {\tilde\eta}=\eta T\leqno (5.18)$$
for function $b\not=0$. It follows that
$$ <d{\tilde y},d{\tilde y}>=\frac{1}{b^2}<dy,dy>.\leqno (5.19)$$
If $\{E_i\}$ is an orthonormal basis for $<dy,dy>$, then
$\{{\tilde E_i}=bE_i\}$ is an orthonormal basis for $<d{\tilde
y},d{\tilde y}>$. From (5.18) and (5.19) we obtain
$$ \sum_i<{\tilde E_i}({\tilde \eta}),{\tilde E_i}({\tilde
\eta})>=b^2\sum_i<E_i(\eta),E_i(\eta)>.\leqno (5.20)$$ It follows
from (5.19) and (5.20) that
$$
g=(\sum_i<E_i(\eta),E_i(\eta)>)<dy,dy>=(\sum_i<E_i(\eta),E_i(\eta)>)III\leqno
(5.21)$$ is a Laguerre invariant. From the last equation of (5.15)
we get
$$ \sum_i<E_i(\eta),E_i(\eta)>=\sum_i(r_i-r)^2.\leqno (5.22)$$
Thus we know that
$$ g=(\sum_i(r_i-r)^2)III\leqno (5.23)$$
is a Laguerre invariant metric at any non-umbilical point of $x$.
We call $g$ the Laguerre metric of $x$. The volume of $g$ is given
by
$$ L(x)=Vol_g(x)=\int_M\frac{(\sum_i(r_i-r)^2)^{(n-1)/2}}{r_1r_2\cdots r_{n-1}}dM,\leqno (5.24)$$
where $dM$ is the volume form with respect to $dx\cdot dx$. We
call critical hypersurfaces of the functional $L(x)$ Laguerre
minimal hypersurfaces. In case $n=3$, we get
$$ L(x)=2\int_M\frac{H^2-K}{K}dM,\leqno (5.25)$$
which is the Laguerre functional given in Blaschke's book [1], the
papers of E.Musso and L.Nicolodi [5] and B. Palmer [8] (up to a
factor).

\newpage\noindent {\bf {\S} 6. Laguerre invariant system for
hypersurfaces in $\R^n$}
\par\medskip
Let $x: M\to \R^n$ be a umbilical free hypersurface with non-zero
principal curvatures. We define
$$ Y=\r\, (\xi,-x\cdot\xi,x\cdot\xi,1),\hskip 5pt\r=\sqrt {\sum_i(r_i-r)^2}>0.\leqno (6.1)$$
If $x,{\tilde x}: M\to \R^n$ are Laguerre equivalent by $T\in
L\G$, then by (5.18), (5.20) and (5.22) we obtain ${\tilde Y}=YT$.
Thus
$$ Y: M\to C^{n+2}\subset \R^{n+3}_2$$
is a Laguerre invariant. We call $Y$ the Laguerre position vector
of the hypersurface $x: M\to \R^n$. It follows immediately from
Theorem 5.1 that
\par\medskip
{\bf Theorem 6.1}{\hskip 3pt}{\it Let $x,{\tilde x}: M^{n-1}\to
\R^n$ be two umbilical free oriented hypersurfaces with non-zero
principal curvatures. Then $x$ and ${ \tilde x}$ are Laguerre
equivalent if and only if there exists $T\in \G$ such that
${\tilde Y}=YT$.}
\par\medskip
Let $Y$ the Laguerre position vector of a hypersurface $x: M\to
\R^n$. Then the Laguerre metric $g$ can be written as
$$ g=<dY,dY>.\leqno (6.2)$$
We denote by $\Delta$ the Laplacian operator of $g$ and define
$$ N=\frac{1}{n-1}\Delta Y+\frac{1}{2(n-1)^2} <\Delta Y,\Delta
Y>Y.\leqno (6.3)$$ Then we have
$$ <Y, Y>=<N,N>=0,\hskip 5pt <Y, N>=-1.\leqno (6.4)$$
Let $\{E_1,E_2,\cdots,E_{n-1}\}$ be an orthonormal basis for
$g=<dY,dY>$ with dual basis $\{\o_1,\o_2,\cdots,\o_{n-1}\}$. Then
we have the following orthogonal decomposition
$$ \R^{n+3}_2=span\{Y,N\}\oplus
span\{E_1(Y),E_2(Y),\cdots,E_{n-1}(Y)\}\oplus \{\eta,\wp\}.\leqno
(6.5)$$ We call $\{Y,N,E_1(Y),E_2(Y),\cdots,E_{n-1}(Y),\eta,\wp\}$
a Laguerre moving frame in $\R^{n+3}_2$ of $x$. By taking
derivatives of this frame we obtain the following structure
equations:
$$ E_i(N)=\sum_jL_{ij}E_j(Y)+C_i\wp\,;\leqno (6.6)$$
$$ E_j(E_i(Y))=L_{ij}Y+\d_{ij}N
+\sum_k\Gamma^k_{ij}E_k(Y)+B_{ij}\wp\,;\leqno (6.7)$$
$$ E_i(\eta)=-C_iY+\sum_jB_{ij}E_j(Y).\leqno (6.8)$$
From these equations we obtain the following basic Laguerre
invariants:
\begin{description}
\item (i) The Laguerre metric $g=<dY,dY>$;
\item (ii) The Laguerre second fundamental form ${\mathbb
B}=\sum_{ij}B_{ij}\o_i\otimes\o_j$;
\item (iii) The symmetric 2-tensor ${\mathbb
L}=\sum_{ij}L_{ij}\o_i\otimes\o_j$;
\item (iv) The Laguerre form $C=\sum_iC_i\o_i$.
\end{description}
By taking further derivatives of (6.6), (6.7) and (6.8) we get the
following relations between these invariants:
$$ L_{ij,k}=L_{ik,j}\,;\leqno (6.9)$$
$$ C_{i,j}-C_{j,i}=\sum_k (B_{ik}L_{kj}-B_{jk}L_{ki})\,;\leqno
(6.10)$$
$$ B_{ij,k}-B_{ik,j}=C_j\d_{ik}-C_k\d_{ij}\,;\leqno (6.11)$$
$$
R_{ijkl}=L_{jk}\d_{il}+L_{il}\d_{jk}-L_{ik}\d_{jl}-L_{jl}\d_{ik}\,;\leqno
(6.12)$$ where $\{L_{ij,k}\}, \{C_{i,j}\}, \{B_{ij,k}\}$ are
covariant derivatives with respect to the Laguerre metric $g$, and
$R_{ijkl}$ is the curvature tensor of $g$. Since $\{E_i'=\r E_i\}$
is an orthonormal basis for the third fundamental form $III$, we
get from (6.8) and (5.22) that
$$ \sum_{ij}B_{ij}^2=\r^{-2}\sum_i<E_i'(\eta),
E_i'(\eta)>=\r^{-2}\sum_i(r_i-r)^2=1.\leqno (6.13)$$ Moreover,
from (6.7) we have
$$ \Delta Y=(\sum_iL_{ii})Y+(n-1)N+(\sum_iB_{ii})\wp.$$
It follows from (6.3) that
$$ \sum_iB_{ii}=0,\hskip 5pt
\sum_i L_{ii}=-\frac{1}{2(n-1)}<\Delta Y,\Delta Y>.\leqno (6.14)$$
From (6.11) we get
$$ \sum_iB_{ij,i}=(n-2)C_j.\leqno (6.15)$$
From (6.12) we get
$$ R_{ik}=-(n-3)L_{ik}-(\sum_iL_{ii})\d_{ik};\leqno (6.16)$$
$$R=-2(n-2)\sum_{i}L_{ii}=\frac{(n-2)}{(n-1)}<\Delta Y,\Delta
Y>.\leqno (6.17)$$
\par\medskip
In case $n>3$, we know from (6.15) and (6.16) that $C_i$ and
$L_{ij}$ are completely determined by the Laguerre invariants
$\{g,\B\}$, thus we get
\par\medskip
{\bf Theorem 6.2}{\hskip 5pt}{\it Two umbilical free oriented
hypersurfaces in $\R^n$ $(n>3)$ with non-zero principal curvatures
are Laguerre equivalent if and only if they have the same Laguerre
metric $g$ and Laguerre second fundamental form $\B$.}
\par\medskip
Let $S: TM\to TM$ be the shape operator for $x$ with principal
radius $r_i$. By direct calculation for hypersurface $x: M\to\R^n$
we get $ B_{ij}=\r^{-1}(r_i-r)\d_{ij}$. We define Laguerre shape
operator
$${\mathbb S}=\r^{-1}(S^{-1}-r\,id): TM\to TM.\leqno (6.18)$$
Then ${\mathbb S}$ is a self-adjoint operator with respect to the
Laguerre metric $g=\r^2\,III$. It follows that for different
principal radius $r_i, r_j$ and $r_k$, the quotient
$(r_i-r_j)/r_i-r_k)$ is a Laguerre invariant.
\par\medskip
In case $n=3$, a complete Laguerre invariant system for surfaces
in $\R^3$ is given by $\{g,\B,{\mathbb L}\}$.
\par\bigskip\noindent
{\bf {\S} 7. Laguerre minimal hypersurfaces in $\R^n$}
\par\medskip
In this section we calculate the first variation formula for
Laguerre minimal hypersurfaces in $\R^n$.
\par\medskip
Let $x_0: M\to\R^n$ be a compact oriented hypersurface in $\R^n$
with boundary $\partial M$. We assume that $x_0$ is umbilical free
and its principal curvatures are non-zero. Let $x: M\times \R\to
\R^n$ be a variation of $x_0$, such that for each $t\in (-\e,\e)$
hypersurface $x_t=x(t,\cdot): M\to \R^n$ is umbilical free and its
principal curvatures are non-zero. Moreover, for any point
$p\in\partial M$ we have
$$ x_t(p)=x_0(p),\hskip 5pt dx_t(p)=dx_0(p): T_pM\to T_pM.\leqno
(7.1)$$ Then the Laguerre volume of $x_t$ is given by
$$ L(t)=L(x_t)=\int_M \frac{(\sum_i(r-r_i)^2)^{(n-1)/2}}{r_1r_2
\cdots r_{n-1}} dM.\leqno (7.2)$$ Our purpose is to calculate the
derivative $L'(0)$.
\par\medskip
Let $\{E_1,E_2,\cdots,E_{n-1}\}$ be an orthonormal basis for the
Laguerre metric $g_t$ of $x_t$ with dual basis
$\{\o_1,\o_2,\cdots,\o_{n-1}\}$ for $T^*M$. We write the variation
vector field of $x: M\times \R\to \R^n$ by
$$ \frac{\partial x}{\partial t}=\sum_iu_ie_i(x)+\r^{-1}f\xi,\leqno
(7.3)$$ where $\xi$ is the unit normal of the hypersurface $x_t$
and $\r$ the function defined by (6.1) for $x_t$. Then we have
$$\frac{\partial \xi}{\partial t}=\sum_i w_iE_i(x)  \leqno (7.4)$$
for some functions $w_i$. Since the second fundamental form is
non-degenerate, we can write $ E_i(x)=\sum_j A_{ij}E_j(\xi)$ with
$det (A_{ij})\not=0$. Let $Y$ be the Laguerre position vector of
$x_t$ given by (6.1). Thus we have
$$ (x\cdot E_i(x),-x\cdot E_i(x),E_i(x),0)=0\hskip 5pt
mod \{Y, E_1(Y),\cdots,E_{n-1}(Y)\}.\leqno (7.5)$$ It follows from
(7.3), (7.4) and (7.5) that
$$ \frac{\partial Y}{\partial t}=\r (\frac{\partial
x}{\partial t}\cdot\xi,-\frac{\partial x}{\partial
t}\cdot\xi,\0,0)=f\wp\,,\hskip 5pt mod \{Y,
E_1(Y),\cdots,E_{n-1}(Y)\}.$$ Thus we can write
$$ \frac{\partial Y}{\partial t}=\s Y+\sum_iv_iE_i(x)+f\wp\leqno
(7.6)$$ for some function $\s$ and some tangent vector field
$V=\sum_iv_iE_i$. We note that the function $f$ is determined by
the normal component of the variation vector field given by (7.3).
\par\medskip
Let $\{Y,N,E_1(Y),\cdots,E_{n-1}(Y),\eta,\wp\}$ be the Laguerre
moving frame of $x_t$. Using the products of the frame we get
$$ dY=\a Y+\sum_i\O_i E_i(Y)+\b \wp\,;\leqno (7.7)$$
$$ dN=-\a N+\sum_i\Psi_iE_i(Y)+\g\wp\,;\leqno (7.8)$$
$$ dE_i(Y)=\Psi_iY+\O_iN+\sum_j\O_{ij}E_j(Y)+\Phi_i\wp\,;\leqno
(7.9)$$
$$ d\eta =-\g Y-\b N+\sum_i\Phi_iE_i(Y),\leqno (7.10)$$
where $\{\a,\b,\O_i,\O_{ij},\Psi_i,\Phi_i,\g)$ are some 1-forms on
$M\times\R$. From (6.6), (6.7), (6.8), (7.6) and the formula
$$ d=\sum_i\o_iE_i(x)+dt\frac{\partial}{\partial t}:
C^{\infty}(M\times\R)\to\Lambda^1(M\times \R)$$
we get
$$ \a=\s dt;\hskip 5pt \O_i=\o_i+v_i dt,\hskip 5pt \b=f dt, \Psi_i=\sum_jL_{ij}\o_j+a_i dt;\leqno (7.11)$$
$$ \Phi_i=\sum_jB_{ij}\o_j+b_i dt; \hskip 5pt
\O_{ij}=\o_{ij}\o_k+p_{ij}dt,\hskip 5pt \g=\sum_iC_i\o_i+c
dt,\leqno (7.12)$$ where $a_i$, $b_i$, $c$, $p_{ij}$ are functions
with $p_{ij}+p_{ji}=0$ and $\o_{ij}$ be the connection form of
$g_t$.
\par\medskip
Taking derivatives of (7.7), (7.8), (7.9) and (7.10) we get
$$ d\b-\sum_i\O_i\wedge\Phi_i-\a\wedge \b=0;\leqno (7.13)$$
$$ d\O_i-\sum_j\O_j\wedge\O_{ji}-\a\wedge\O_i=0;\leqno (7.14)$$
$$
d\Phi_i-\sum_j\O_{ij}\wedge\Phi_j-\Psi_i\wedge\b-\O_i\wedge\g=0.\leqno
(7.15)$$ Since
$$ d=d_M+dt\wedge\frac{\partial}{\partial t}:
\Lambda^1(M\times\R)\to\Lambda^2(M\times \R),$$ where $d_M$ is the
differential operator on $M$, by comparing the coefficients of
$\o_i\wedge dt$ of (7.13) and (7.14) we get
$$ b_i=E_i(f)+\sum_jB_{ij}v_j;\leqno (7.16)$$
$$ \frac{\partial\o_i}{\partial
t}=\sum_j(v_{i,j}+p_{ij}+\s\d_{ij})\o_j;\leqno (7.17)$$ where
$\{v_{i,j}\}$ is the covariant derivative of $V=\sum_iv_iE_i$. By
comparing the coefficients of $\o_i\wedge dt$ of (7.15) and using
(7.16), (7.17) we get
$$  \frac{\partial B_{ij}}{\partial t}+\s
B_{ij}=f_{i,j}+\sum_kv_k(B_{ki,j}+C_j\d_{ik})
+\sum_k(p_{ik}B_{kj}-p_{jk}B_{ki})-L_{ij}f-c\d_{ij},\leqno
(7.18)$$ where $(f_{i,j})$ is the Hessian matrix of $f$. From
(6.11), (6.13) and (6.14) we have
$$ \sum_iB_{ii}=0,\hskip 5pt \sum_{ij}B_{ij}^2=1,\hskip 5pt
\sum_{ij}(B_{ki,j}+C_j\d_{ik})B_{ij}=0.\leqno (7.19)$$ By
multiplying $B_{ij}$ to (7.18), taking sum and using (7.19) we get
$$ \s=\sum_{ij}(f_{i,j}B_{ij}-fL_{ij}B_{ij}).\leqno (7.20)$$
\par\medskip
Now we come to calculate $L'(0)$. Since the Laguerre volume can be
written as
$$ L(t)=\int_M\o_1\wedge\o_2\wedge\cdots\wedge\o_{n-1},$$
we get from (7.17) that
$$ L'(0)=\int_M (div V+(n-1)\s)dM,$$
where $V=\sum_iv_iE_i$. Since on $\partial M$ we have $v_i=0$,
$f=0$ and $f_i=E_i(f)=0$, it follows from (7.20) and Green-formula
that
$$ L'(0)=(n-1)\int_M\sum_{ij}(B_{ij,ij}-L_{ij}B_{ij})f
dM.\leqno (7.21)$$ Thus we obtain
\par\medskip
{\bf Theorem 7.1}{\hskip 3pt}{\it The first variation of a
Laguerre volume for hypersurfaces in $\R^n$ depends only on the
normal component of the variation vector field. The Euler-Lagrange
equation of the Laguerre functional is given by}
$$ \sum_{ij}(B_{ij,ij}-L_{ij}B_{ij})=0.\leqno (7.22)$$
\par\medskip
Using (6.15) we can write the Euler-Lagrange equation by
$$ \sum_iC_{i,i}-\frac{1}{n-2}\sum_{ij}L_{ij}B_{ij}=0.\leqno
(7.23)$$ From (6.8), (6.7) and (6.15) we obtain
$$\Delta \eta= \sum_i (-C_{i,i}+\sum_jL_{ij}B_{ij})Y+(n-3)\sum_i
C_iE_i(Y)+\wp.\leqno (7.24)$$
\par\bigskip\noindent
{\bf {\S} 8. Hypersurfaces in the Laguerre space forms}
\par\medskip
Using Laguerre embedding $\s$ and $\tau$ in $\S 4$ we can regard a
hypersurface $x: M\to\R^n_1$ or $x: M\to\R^n_0$ as a Laguerre
hypersurface $(x',\xi'): M\to U\R^n$. In this section we study the
relations between $x$ and $x'$.
\par\medskip
Let $x: M\to \R^n_1$ be a space-like oriented hypersurfaces in the
Lorentzian space $\R^n_1$ with the inner product $<\,,\,>$ given
in (4.1). Let $\xi$ be the normal of $x$ with $<\xi,\xi>=-1$. The
shape operator $S: TM\to TM$ of $x$ is defined by $ d\xi =-dx\circ
S$. Since $S$ is self-adjoint on $TM$, all eigenvalues $\{k_i\}$
of $S$ are real. We assume that $k_i\not=0$. We define $r_i=1/k_i$
the curvature radius of $x$ and $r=(r_1+r_2+\cdots+r_{n-1})/(n-1)$
the mean curvature radius of $x$. Let $e_i$ be the eigenvector of
$x$ with respect to the eigenvalue $k_i$. Then we have
$$ e_i(x)=-r_ie_i(\xi).\leqno (8.1)$$
In particular, we have $e_i(x_1)=-r_ie_i(\xi_1)$. We define $
(x',\xi')=\s (x,\xi): M\to U\R^n$, where $\s: U\R^n_1\to U\R^n$ is
the Laguerre embedding given by (4.11). By a direct calculation we
get from (4.11) and (8.1) that
$$ e_i(x')=-(r_i\xi_1+x_1)e_i(\xi').\leqno (8.2)$$
It follows that $e_i$ is the principal vector for the Laguerre
hypersurface $f'=(x',\xi'): M\to U\R^n$ corresponding to the
curvature radius
$$ r_i'=r_i\xi_1+x_1, \leqno (8.3)$$
which implies the following relations between the mean curvature
radius $r'$ and $r$:
$$ r'=r\xi_1+x_1,\hskip 5pt  \r\,'^2=\sum_i(r'_i-r')^2=\xi_1^2\sum_i(r_i-r)^2=\xi_1^2\r^2.\leqno (8.4)$$
It is easy to verify from (4.11) that
$$ Y'=\r'(x'\cdot\xi',-x'\cdot\xi',\xi',1)=\r(<x,\xi>,-<x,\xi>,1,\xi>)=Y.\leqno (8.5)$$
Thus the Laguerre metric is given by
$$ g'=\r'\,^2III'=<dY',dY'>=<dY,dY>=\r^2III=g,\leqno (8.6)$$
where $III$ and $III'$ are the third fundamental forms of $x$ and
$x'$, respectively. By (4.4) we know that the mean curvature
radius sphere $H(x+r\xi,r)$ in $U\R^n_1$ corresponds to the vector
$[\eta]\in\Q^{n+1}$, where
$$ \eta =(\frac{1}{2}(1+<x,x>),\frac{1}{2}(1-<x,x>),0,
x)+r\,(<x,\xi>, -<x,\xi>,1,\xi).\leqno (8.8)$$ By a direct
calculation we know that
$$ \eta' =(\frac{1}{2}(1+|x'|^2),\frac{1}{2}(1-|x'|^2),
x',0)+r'\,(x'\cdot\xi', -x'\cdot\xi',\xi',1)=\eta .\leqno (8.9)$$
Thus the Laguerre embedding $\s: U\R^n_1\to U\R^n$ takes the mean
curvature radius sphere $H(x+r\xi,r)$ in $\R^n_1$ into the mean
curvature radius sphere $S(x'+r'\xi',r')$ in $\R^n$.
\par\medskip
Let $x: M\to \R^n_0$ be a space-like oriented hypersurfaces in the
degenerate hyperplane $\R^n_0$ of the Lorentzian space
$\R^{n+1}_1$ with the inner product $<\,,\,>$ given in (4.12). Let
$\xi$ be the unique vector satisfying
$$ <\xi, dx>=0,\hskip 5pt <\xi,\xi>=0,\hskip 5pt <\xi,\nu>=1.\leqno (8.10)$$
We define the shape operator $S: TM\to TM$ by $d\xi=-dx\circ S$.
Since $S$ is self-adjoint, all eigenvalues $\{k_i\}$ of $S$ are
real. We assume that $k_i\not=0$. Then we define $r_i=1/k_i$ the
curvature radius of $x$ and $r=(r_1+r_2+\cdots+r_{n-1})/(n-1)$ the
mean curvature radius of $x$. Let $\tau: U\R^n_0\to U\R^n$ be the
Laguerre embedding defined by (4.23) and $(x',\xi')=\tau\circ
(x,\xi)$. By a similar way as we can show that
$$ Y=\r\,(<x,\xi>,-<x,\xi>,\xi)=\r'\,(x'\cdot\xi',
-x'\cdot\xi',\xi',1)=Y',\leqno (8.11)$$ where
$\r'\,^2=\sum_i(r_i'-r')^2$ and $\r^2=\sum_i(r_i-r)^2$. Thus the
Laguerre metric is given by
$$ g=\r^2III=\r'\,^2III'=g', \leqno (8.12)$$
where $III$ and $III'$ are the third fundamental forms of $x$ and
$x'$, respectively. Moreover, we have $\eta=\eta'$, where
$$ \eta=(\frac{1}{2}(1+<x,x>),\frac{1}{2}(1-<x,x>),x)+r
(<x,\xi>,-<x,\xi>,\xi);\leqno (8.13)$$
$$ \eta' =(\frac{1}{2}(1+|x'|^2),\frac{1}{2}(1-|x'|^2),
x',0)+r'\,(x'\cdot\xi', -x'\cdot\xi',\xi',1).\leqno (8.14)$$
\par\medskip
It follows immediately from (6.1), (5.31), (8.5), (8.8), (8.11)
and (8.13) that
\par\medskip
{\bf Proposition 8.1}{\hskip 5pt}{\it Let $x$ be a Laguerre
hypersurface in $U\R^n$, $U\R^n_1$ or $U\R^n_0$. Let $\c\in
\R^{n+3}_2$ be the time-like vector $\c=(0,0,\0,-1)$, the
space-like vector $\c=(0,0,1,\0)$ or the light-like
$\c=(0,0,\nu)$, respectively. Let $\{r_i\}$ be the curvature
radius of $x$. Let $r$  the mean curvature radius of $x$ and
$\r^2=\sum_i(r_i-r)^2$. Then we have $<Y,\c>=\r$ and
$<\eta,\c>=r$.}
\par\medskip
{\bf Theorem 8.1}{\hskip 5pt}{\it A surface in $\R^3$, $\R^3_1$ or
$\R^3_0$ (regarded as a Laguerre surface in $\R^3$) is Laguerre
minimal if and only if $\Delta_{III}\,r=0$. }
\par\noindent
{\bf Proof.}{\hskip 5pt} From (7.24) we know that for surface case
we have
$$ \Delta_{III} \eta=\r^2 \sum_i (-C_{i,i}+\sum_jL_{ij}B_{ij})Y+\r^2\wp.\leqno (8.15)$$
Any surface in $\R^3_1$ and $\R^3_0$ can be regarded as surface in
$\R^3$ with the same $Y$ and $\eta$. Let $\c\in \R^{n+3}_2$ be the
vector given in Proposition 8.1. Since $<\c,\wp>=0$, we get from
(8.15) that
$$ \Delta_{III} r= \r^3 \sum_i (-C_{i,i}+\sum_jL_{ij}B_{ij}).\leqno (8.16)$$
Thus we complete the proof of Theorem 8.1.
\par\medskip
{\bf Remark.}{\hskip 5pt}{\it For Laguerre minimal surface in
$\R^3$ the equation $\Delta_{III} r=0$ is given in Blaschke's book
[1].}
\par\medskip
Since $r=(k_1+k_2)/k_1k_2$, we know that $r=0$ if and only if $x$
is a minimal surface in $\R^3$, $\R^3_1$ or $\R^3_0$. Thus a
minimal surface in $\R^3_1$ or $\R^3_0$ induces a Laguerre minimal
surface in $U\R^3$ (by using Laguerre embedding).
\par\medskip
{\bf Theorem 8.2}{\hskip 5pt}{\it The only compact Laguerre
minimal surface in $\R^3$ is the round sphere.}
\par\noindent
{\bf Proof.}{\hskip 5pt} Let $x: M\to\R^3$ be a compact Laguerre
minimal surface. From (8.15) we $\Delta_{III}\eta =\r^2\wp$, which
holds also at umbilical points of $x$. Since $\eta
=(\a,\b,x+r\xi,r)$ for some functions $\a$ and $\b$, we get
$$ \Delta_{III}r=0,\hskip 5pt \Delta_{III}(x+r\xi)=0.$$
Since $M$ is compact, we know that both $r$ and $x+r\xi=x_0$ are
constant. Thus we have $|x-x_0|^2=r^2$ and $x$ is a round sphere
in $\R^3$. We complete the proof of Theorem 8.2.

\begin{flushleft}
Department of Mathematics \\
Beijing Institute of Technology\\
Beijing 100081, P. R. China\\
\texttt{litz@163.com}\\
\vspace{5mm}
LMAM\\
School of Mathematical Sciences\\
Peking University\\
Beijing 100871,P. R. China\\
\texttt{cpwang@math.pku.edu.cn}
\end{flushleft}

\end{document}